\documentclass[10pt,reqno,fleqn]{amsart}
\usepackage{amsmath,amscd,amssymb,amsfonts,amsthm,courier,relsize,bm}
\usepackage{hyperref,enumitem,mathrsfs,slashed,mathtools, extarrows}
\usepackage[bottom]{footmisc}
\usepackage{xcolor} 
\usepackage[all]{xy}
\usepackage{a4wide}
\usepackage[bitstream-charter]{mathdesign} %In this one, replace \nu with \!\nu
\usepackage[T1]{fontenc}
\usepackage{tikz}
\usepackage[bottom]{footmisc}
\usepackage{hyperref}

\hypersetup{colorlinks,linkcolor={blue},citecolor={blue},urlcolor={red}}  

\numberwithin{equation}{section}
\newtheoremstyle{mytheoremstyle} % name
    {4mm}                    % Space above
    {4mm}                    % Space below
    {\itshape}                   % Body font
    {6mm}                           % Indent amount
    {\scshape}                   % Theorem head font
    {.}                          % Punctuation after theorem head
    {0.5em}                       % Space after theorem head
    {}  % Theorem head spec (can be left empty, meaning ?normal?)

\theoremstyle{mytheoremstyle}

\newtheorem{df}{Definition}[section]
\let\olddf\df
\renewcommand{\df}{\olddf\normalfont}
\newtheorem{thm}[df]{Theorem}
\newtheorem{prop}[df]{Proposition}

\newtheorem{lem}[df]{Lemma}
\newtheorem{cor}[df]{Corollary}
\newtheorem{ex}[df]{Example}
\let\oldex\ex
\renewcommand{\ex}{\oldex\normalfont}
\newtheorem{rk}[df]{Remark}
\let\oldrk\rk
\renewcommand{\rk}{\oldrk\normalfont}
\newtheorem*{pr}{Proof}
\let\oldpr\pr
\renewcommand{\pr}{\oldpr\normalfont}

\newcommand{\C}{\mathbb{C}}
\newcommand{\R}{\mathbb{R}}
\newcommand{\Z}{\mathbb{Z}}
\newcommand{\N}{\mathbb{N}}

\newcommand{\A}{\mathcal{A}}

\renewcommand{\S}{S}
\newcommand{\D}{\mathcal{D}}
\newcommand{\B}{\mathcal{B}}
\newcommand{\E}{\mathcal{E}}
\newcommand{\Tr}{\mathrm{Tr}}
\newcommand{\cinf}{C^\infty}

\newcommand{\Ind}{\mathrm{Ind}}

\newcommand{\ch}{\mathrm{ch}}

\newcommand{\ad}{\mathrm{ad}}

\renewcommand{\dim}{\mathrm{dim}}

\renewcommand{\Re}{\mathrm{Re}}

\newcommand{\Res}{\mathrm{Res}}

\newcommand{\dom}{\mathrm{dom}}

\renewcommand{\ch}{\mathrm{ch}}

\newcommand{\HP}{\mathrm{HP}}
\newcommand{\Pf}{\mathrm{Pf}}

\newcommand*{\barint}{\mathop{\ooalign{$\displaystyle{\int}$\cr$-$}}}
\newcommand*{\littlebarint}{\mathop{\ooalign{$\int$\cr$-$}}}
\renewcommand{\dom}{\mathrm{dom}}

\newcommand{\Op}{\mathrm{Op}}
\newcommand{\alg}{\mathrm{alg}}
\renewcommand{\\}{\vspace{2mm}}
\renewcommand{\i}{\textup{\textbf{i}}}

\newcommand{\Diff}{\mathrm{Diff}}
\renewcommand{\tilde}{\widetilde}
\renewcommand{\epsilon}{\varepsilon}
\renewcommand{\sim}{\thicksim}

\title{The Radul cocycle, the Chern--Connes character, and manifolds with conical singularities}
\author{Rudy Rodsphon}
\address{Department of Mathematics and Statistics \\ Washington University in St. Louis \\  1 Brookings Dr., Cupples I Hall \\ St Louis, MO 63130}
\email{rrudy@wustl.edu}

\begin{document}

\usetikzlibrary{arrows,chains,matrix,positioning,scopes}
\makeatletter
\tikzset{join/.code=\tikzset{after node path={%
\ifx\tikzchainprevious\pgfutil@empty\else(\tikzchainprevious)%
edge[every join]#1(\tikzchaincurrent)\fi}}}
\makeatother
\tikzset{>=stealth',every on chain/.append style={join},
         every join/.style={->}}
\tikzstyle{labeled}=[execute at begin node=$\scriptstyle,
   execute at end node=$]

\definecolor{c3}{RGB}{0, 255, 180}
\definecolor{c2}{RGB}{0, 155, 255}
\definecolor{c1}{RGB}{255, 13, 214}
\definecolor{c23}{RGB}{0,205,218}
\definecolor{c123}{RGB}{128,109,216}

\definecolor{BleuICJ}{RGB}{105,166,207}

%\begin{abstract}

%\setlength{\parindent}{0mm}
%\textsc{Keywords.} Cyclic cohomology, K-theory, Index theory, Pseudodifferential operators \\

%\textsc{MSC.}  19D55, 19K56, 58J42, 46L87
%\end{abstract}

\maketitle 

\setlength{\parindent}{6mm}

\section*{Introduction}

In a former article \cite{Rod2013}, we had established a formula calculating the Chern character (in K-homology) of an abstract pseudodifferential extension in terms of residues of zeta functions, applicable in the presence of multiple poles. We direct the reader to the aforementioned article for further references on the history of such a formula (see e.g \cite{Nis1997}), and recent applications (see e.g \cite{Per2012, Per2014, PerRod2014}). After establishing a precise relationship between this formula and the residue cocycle of Connes--Moscovici \cite{CM1995}, we discuss briefly the application of this formula in the context of manifolds with conical singularities, which may exhibit triple poles. Let us give an overview of this note.  \\

We first recall some material from \cite{Hig2006}, which develops a formalism of algebras of abstract differential operators, and provides a relationship with regular spectral triples. This yields naturally an abstract pseudodifferential extension, and we then recall the derivation of its Chern character in the form of a cyclic 1-cocycle generalizing the Radul cocycle, applicable contexts where the zeta function exhibits multiple poles. We then discuss the application of this formula to the context of manifolds with conical singularities, and the associated spectral triples.

\setlength{\parindent}{0mm}

\section{Abstract Differential Operators and Traces} \label{Abstract Differential Operators and a Trace}

In this part, we recall the Abstract Differential Operators formalism developed by Higson in \cite{Hig2006} to simplify the proof of the Connes-Moscovici local index formula \cite{CM1995}. For more details, one can refer to Higson's article or \cite{Rod2013}.   

\subsection{Abstract Differential Operators} \label{Abstract Differential Operators}

Let $H$ be a (complex) Hilbert space and let $\Delta$ be an unbounded, positive and self-adjoint operator acting on it, with domain $\dom{\Delta}$. To keep the exposition simple, we suppose that $\Delta$ has a compact resolvent. \\

We denote by $H^\infty$ the following intersection:
\begin{equation*}  
H^\infty = \bigcap_{k=0}^{\infty} \dom(\Delta^k)
\end{equation*}

\begin{df} \label{df ADO}
An algebra $\mathcal{D}(\Delta)$ of \emph{abstract differential operators} associated to $\Delta$ is an algebra of operators on $H^\infty$ satisfying the following conditions
\begin{enumerate}[label={(\roman*)}]
\item The algebra $\D(\Delta)$ is filtered,
\begin{equation*} 
\D(\Delta) = \bigcup_{q=0}^{\infty} \D_q(\Delta)
\end{equation*}
\end{enumerate}
that is $\D_p(\Delta) \cdot \D_q(\Delta) \subset \D_{p+q}(\Delta)$. An element $X \in \D_q(\Delta)$ is an \emph{abstract differential operator of order at most $q$}.  \\

\begin{enumerate}[label={(\roman*)}]
\item[(ii)] There is a $r>0$ ("the order of $\Delta$") such that for every $X \in \D_q(\Delta)$, $[\Delta,X] \in \D_{r+q-1}(\Delta)$. \\
\end{enumerate}

For $s \in \R$, define the $s$-\emph{Sobolev space} $H^s$ as the subspace $\dom(\Delta^{s/r})$ of $H$, which is a Hilbert space when equipped with the norm 
\[ \Vert v \Vert_s   = (\Vert v \Vert^2 + \Vert \Delta^{s/r}v \Vert^2)^{1/2} \]  

\begin{enumerate}[label={(\roman*)}]
\item[(iii)] \emph{Elliptic estimate.} If $X \in \D_{q}(\Delta)$, then, there is a constant $\varepsilon > 0$ such that 
\begin{equation*} 
\Vert v \Vert_q + \Vert v \Vert \geq \varepsilon \Vert Xv \Vert \, , \, \forall v \in H^\infty 
\end{equation*}
\end{enumerate}
\end{df}

Having G\"{a}rding's inequality in mind, the elliptic estimate exactly says that $\Delta^{1/r}$ should be thought as an "abstract elliptic operator" of order 1. It also says that any differential operator $X$ of order $q$ can be extended to a bounded operator form $H^{s+q}$ to $H^{s}$. This last property will be useful to define pseudodifferential calculus in this setting. The main example to keep in mind is of course the case in which $\Delta$ is a Laplace type operator on a closed Riemannian manifold $M$. 

\subsection{Correspondence with spectral triples}

Let $(A,H,D)$ a spectral triple (cf. \cite{CM1995} or \cite{Hig2006}). One may construct a algebra of abstract differential operators $\D = \D(A,D)$ recursively as follows : 
\begin{gather*} 
\D_0 = \text{algebra generated by } A \text{ and } [D,A] \\
\D_1 = [\Delta, \D_0] + \D_0[\Delta, \D_0] \\
\qquad \vdots \\
\D_k = \sum_{j=1}^{k-1} \D_j \cdot \D_{k-j} + [\Delta, \D_{k-1}] + \D_0[\Delta, \D_{k-1}]
\end{gather*}

Let $\delta$ be the unbounded derivation $\ad \vert D \vert = [\vert D \vert, \, . \,]$ on $\B(H)$. The spectral triple is $(A,H,D)$ is said \emph{regular} if $A, [D,A]$ are included in $\bigcap_{n=1}^{\infty} \dom \, \delta^{n}$. The following theorem of Higson relates algebras of abstract differential operators and spectral triples.

\begin{thm} \emph{(Higson, \cite{Hig2006}).} Suppose that $A$ maps $H^\infty$ into itself. Then, the spectral triple $(A,H,D)$ is regular if and only if the elliptic estimate of Definition \ref{df ADO} holds. 
\end{thm}

\subsection{Zeta Functions} \label{Zeta Functions}

Let $\D(\Delta)$ be an algebra of abstract differential operators. For $z \in \C$, one defines the \emph{complex powers} $\Delta^{-z}$ of $\Delta$ using functional calculus :
\begin{equation*}
\Delta^{-z} = \frac{1}{2 \pi i} \int \lambda^{-z} (\lambda - \Delta)^{-1} d\lambda
\end{equation*}
where the contour of integration is a vertical line pointing downwards separating $0$ and the (discrete) spectrum of $\Delta$. This converges in the operator norm when $\Re(z) > 0$, and using the semi-group property, all the complex powers can be defined after multiplying by $\Delta^k$, for $k \in \N$ large enough. Moreover, since $\Delta$ has compact resolvent, the complex powers of $\Delta$ are well defined operators on $H^\infty$. \\

We will suppose that there exists a $d\geq 0$ such that for every $X \in \D_q(\Delta)$, the operator $X \Delta^{-z}$ extends to a trace-class operator on $H$ for $z$ on the half-plane $\Re(z) > \frac{q+d}{r}$. The \emph{zeta function} of $X$ is 
\[ \zeta_X(z) = \Tr(X \Delta^{-z/r}) \]
The smallest $d$ verifying the above property is called the \emph{analytic dimension} of $\D(\Delta)$. In this case, the zeta function is holomorphic on the half-plane $\Re(z) > q+d$. We shall say that $\D(\Delta)$ has the \emph{analytic continuation property} if for every $X \in \D(\Delta)$, the associated zeta function extends to a meromorphic function of the whole complex plane. \\

There properties are set for all the section, unless if it is explicitly mentioned. \\

These notions come from properties of the zeta function on a closed Riemannian manifold $M$ : it is well-known that the algebra of differential operators on $M$ has analytic dimension $\dim \,M$ and the analytic continuation property. Its extension to a meromorphic function has at most simple poles at the integers smaller that $\dim \,M$. In the case where $M$ is foliated, the dimension of the leaves appears in the analytic dimension when working in the suitable context. Hence, the zeta function provide informations not only on the topology of $M$, but also on its the geometric structure,  illustrating the relevance of this abstraction.   

\subsection{Abstract Pseudodifferential Operators}  \label{Abstract Pseudodifferential Operators}

Let $\D(\Delta)$ an algebra of abstract differential operators of analytic dimension $d$. To define the notion of pseudodifferential operators, we need a more general notion of order, not necessary integral, which covers the one induced by the filtration of $\D(\Delta)$.

\begin{df}
An operator $T : H^\infty \to H^\infty$ is said to have \emph{pseudodifferential order} $m \in \R$ if for every $s \geq 0$, it extends to a bounded operator from $H^{m+s}$ to $H^s$. In addition, we require that operators of analytic order stricly less than $-d$ are trace-class operators.   
\end{df}

That this notion of order covers the differential order is due to the elliptic estimate, as already remarked in Section \ref{Abstract Differential Operators}. The space of such operators, denoted $\Op(\Delta)$, forms a $\R$-filtered algebra. There is also a notion of regularizing operators which are, as expected, the elements of the (two-sided) ideal of operators of all order. 

\begin{rk} Higson uses in \cite{Hig2006} the term "analytic order", but as the examples we deal with in the paper are about pseudodifferential operators, we prefer the term pseudodifferential order. %Usually, the last condition of the previous definition is not required. Operators between classical Sobolev spaces verify this property, and it will also be the case in the examples we consider through the paper. As we need to consider traces of regularizing operators, we cannot get rid of this condition. 
\end{rk} 

\begin{ex} For every $\lambda \in \C$ not contained in the spectrum of $\Delta$, the resolvent $(\lambda - \Delta)^{-1}$ has analytic order $r$. Moreover, by spectral theory, its norm as an operator between Sobolev spaces is a $O(\vert \lambda \vert^{-1})$. 
\end{ex} 

The following notion is due to Uuye, cf. \cite{Uuy2009}. We just added an assumption on the zeta function which is necessary for what we do.  

\begin{df} \label{abstract pdo} An algebra of abstract pseudodifferential operators is a $\R$-filtered subalgebra $\Psi(\Delta)$ of $\Op(\Delta)$, also denoted $\Psi$ when the context is clear, satisfying 
\begin{align*}
& \Delta^{z/r} \Psi^m \subset \Psi^{\Re(z) + m},  \quad \Psi^m \Delta^{z/r} \subset \Psi^{\Re(z) + m}
\end{align*}
and which commutes, up to operators of lower order, with the complex powers of $\Delta^{1/r}$, that is , for all $m \in \R$, $z \in \C$
\begin{equation*} [\Delta^{z/r}, \Psi^m] \subset \Psi^{\Re(z) + m - 1} \end{equation*}

Moreover, we suppose that for every $P \in \Psi^m(\Delta)$, the zeta function 
\begin{equation*} 
\zeta_P(z) = \Tr(P \Delta^{-z/r}) 
\end{equation*}
is holomorphic on the half-plane $\Re(z) > m+d$, and extends to a meromorphic function of the whole complex plane. %The poles, which are at most of the same order as those of the zeta function, are located in the set 
%\[ \{m+d, m+d-1, \ldots \} \]
We shall denote by 
\begin{equation*} \Psi^{-\infty} = \bigcap_{m \in \R} \Psi^m \end{equation*}  
\end{df}

%As a consequence, we have the following proposition that we shall not prove here. For this, it suffices to adapt the proof of the proposition 4.14 in \cite{Hig2006}. 

Of course, this is true for the algebra of (classical) pseudodifferential operators on a closed manifold. We shall recall later what happens in the example of Heisenberg pseudodifferential calculus on a foliation, as described by Connes and Moscovici in \cite{CM1995}.

We end this part with a notion of asymptotic expansion for abstract pseudodifferential operators. This can be seen as "convergence under the residue".

\begin{df} \label{asymptotic expansion} Let $T$ and $T_\alpha$ ($\alpha$ in a set $A$) be operators on $\Psi$. We shall write 
\begin{equation*} T \sim \sum_{\alpha \in A} T_\alpha \end{equation*} 
if there exists $c > 0$ and a finite subset $F \subset A$ such that for all finite set $F' \subset A$ containing $F$, the map
\[ z \longmapsto \Tr\left( (T - \sum_{\alpha \in F'} T_\alpha) \Delta^{z/r} \right)\]
is holomorphic in a half-plane $\Re(z) > -c$ (which contains $z=0$).
\end{df}

\begin{ex} Suppose that that for every $M > 0$, there exists a finite subset $F \subset A$ such that 
\begin{equation*} T - \sum_{\alpha \in F} T_\alpha \in \Psi^{-M} \end{equation*}
Then, $ T \sim \sum_{\alpha \in A} T_\alpha $
\end{ex}

In this context, asymptotic means that when taking values under the residue, such infinite sums, which have no reason to converge in the operator norm, are in fact finite sums. To this effect, the following lemma is crucial. \\

\begin{lem} \label{CM trick} \emph{(Connes-Moscovici's trick, \cite{CM1995, Hig2006})} Let $Q \in \Psi(\Delta)$ be an abstract pseudodifferential operator. Then, for any $z \in \C$, we have 
\begin{equation}
[\Delta^{-z}, Q] \sim \sum_{k\geq 1} \binom{-z}{k} Q^{(k)}\Delta^{-z-k}
\end{equation}
where we denote $Q^{(k)} = \ad(\Delta)^k(Q)$, $\ad(\Delta) = [\Delta, \, . \,]$.
\end{lem}

The proof relies of the following identity, for $z \in \C$ with $\Re(z)$ large enough (cf \cite{Hig2006}, Lemma 4.20) :
\begin{equation}  \label{CM trick bis}
\Delta^{-z} Q - Q \Delta^{-z} = \sum_{k = 1}^N \binom{-z}{k} Q^{(k)}\Delta^{-z-k} + \frac{1}{2 \pi i} \int \lambda^{-z}(\lambda - \Delta)^{-1} Q^{(N+1)}(\lambda - \Delta)^{-N-1} \, d\lambda
\end{equation}

\subsection{Higher traces on the algebra of abstract pseudodifferential operators} \label{Residue trace}

We give in this paragraph a simple generalization of the Wodzicki residue trace in the case where the zeta function of the algebra $\D(\Delta)$ has poles of arbitrary order. This was already present in the work of Connes and Moscovici (see \cite{CM1995}). 

\begin{prop} \label{higher WG trace}
Let $\Psi(\Delta)$ an algebra of abstract pseudodifferential operators, following the context of the previous paragraphs. Suppose that the associated zeta function has a pole of order $p \geq 1$ in $0$. Then, the functional 
\begin{equation*} \barint^{p} P = \Res_{z=0} z^{p-1}\Tr(P \Delta^{-z/r})  \end{equation*}
defines a trace on $\Psi(\Delta)$.
\end{prop}

\section{The Radul cocycle for abstract pseudodifferential operators}

\subsection{Abstract index theorems} \label{Abstract index}

Let $A$ be an associative algebra over $\C$, possibly without unit, and $I$ an ideal in $A$. The extension
\begin{equation*} 0 \to I \to A \to A /I \to 0 \end{equation*}
gives rise to the following diagram relating algebraic K-theory and periodic cyclic homology 
\begin{equation} \label{Abstract index diagram}
\xymatrix{
    K_1^{\alg}(A /I) \ar[r]^{\Ind} \ar[d]^{\ch_1} & K_0^{\alg}(I) \ar[d]^{\ch_0} \\
    \HP_1(A /I) \ar[r]^{\partial} & \HP_0(I)}
\end{equation}
The vertical arrows are respectively the odd and even Chern character. \\
 
We still denote $\partial :  \HP^0(I) \rightarrow \HP^1(A /I)$ the boundary map in cohomology. As mentioned in \cite{Nis1997}, for $[\tau] \in \HP^0(I)$, $[u] \in K^1(A/I)$, one has the equality :
\begin{equation} \label{index pairing}
\langle [\tau], \ch_0 \Ind [u] \rangle = \langle \partial [\tau], \ch_1 [u] \rangle 
\end{equation}

A standard procedure to calculate the boundary map $\partial$ in cohomology associated to the extension as follows. If $[\tau] \in \HP^0(I)$ is represented by a hypertrace $\tau : I \to \C$, i.e a linear map satisfying the condition $\tau([A,I]) = 0$, then choose a lift $\tilde{\tau} : A \to \C$ of $\tau$, such that $\tilde{\tau}$ is linear (in general, this is not a trace), and a linear section $\sigma : A/I \to A$ such that $\sigma(1) = 1$, after adjoining a unit where we have to. Then, $\partial[\tau]$ is represented by the following cyclic 1-cocycle :
\begin{equation*} 
c(a_0,a_1) = b\tilde{\tau}(\sigma(a_0), \sigma(a_1)) = \tilde{\tau}([\sigma(a_0), \sigma(a_1)])
\end{equation*} 
where $b$ is the Hochschild coboundary.

\subsection{The generalized Radul cocycle}

We can finally come to the main theorem of this section. Let $\D(\Delta)$ be an algebra of abstract differential operators and $\Psi$ be an algebra of abstract pseudodifferential operators. We consider the extension 
\begin{equation*} 0 \to \Psi^{-\infty} \to \Psi \to \S \to 0 \end{equation*}
where $\S$ is the quotient $\Psi / \Psi^{-\infty} $.
The operator trace on $\Psi^{-\infty}$ is well defined, and $\Tr([\Psi^{-\infty}, \Psi]) = 0$. 

\begin{thm} \label{local index formula} Suppose that the pole in zero of the zeta function is of order $p \geq 1$. Then, the cyclic 1-cocycle $\partial [\Tr] \in \HP^1(S)$ is represented by the following functional :
\begin{equation*} 
c(a_0,a_1) = \barint^1 a_0 \delta(a_1) - \dfrac{1}{2!}\barint^2 a_0 \delta^2(a_1) + \ldots + \dfrac{(-1)^{p-1}}{p!}\barint^p a_0 \delta^p(a_1) 
\end{equation*}
where $\delta(a) = [\log \Delta^{1/r}, a]$ and $\delta^k(a) = \delta^{k-1}(\delta(a))$ is defined by induction. We shall call this cocycle as the \emph{generalized Radul cocycle}. 
\end{thm} 

Here, the commutator $[\log \Delta^{1/r}, a]$ is defined as the non-convergent asymptotic expansion 
\begin{equation} \label{log commutator 1}
[\log \Delta^{1/r}, a] \sim \dfrac{1}{r} \sum_{k\geq 1} \frac{(-1)^{k-1}}{k} a^{(k)} \Delta^{-k} 
\end{equation}
where $a^{(k)}$ has the same meaning as in Lemma \ref{CM trick}. 
This expansion arises by first using functional calculus :
\[ \log \Delta^{1/r} = \dfrac{1}{2 \pi \i} \int \log \lambda^{1/r} (\lambda - \Delta)^{-1} \, d\lambda \]
and then, reproducing the same calculations made in the proof of Lemma \ref{CM trick} to obtain the formula (cf. \cite{Hig2006} for details). In particular, note that $\log \Delta^{1/r} = \frac{1}{r} \log \Delta$.  \\

Another equivalent expansion possible, that we will also use, is the following
\begin{equation} \label{log commutator 2}
[\log \Delta^{1/r}, a] \sim \sum_{k\geq 1} \frac{(-1)^{k-1}}{k} a^{[k]} \Delta^{-k/r} 
\end{equation}
where $a^{[1]} = [\Delta^{1/r}, a]$, and $a^{[k+1]} = [\Delta^{1/r}, a^{[k]}]$. 
A heuristic explanation is the following. We first lift the trace on $\Psi^{-\infty}$ to a linear map $\tilde{\tau}$ on $\Psi$ using a zeta function regularization by "Partie Finie" :
\begin{equation*} \tilde{\tau}(P) = \Pf_{z=0} \Tr(P \Delta^{-z/r}) \end{equation*}
for any $P \in \Psi$. The "Partie Finie" $\Pf$ is defined as the constant term in the Laurent expansion of a meromorphic function. Let $Q \in \Psi$ be another pseudodifferential operator. Then, we have
\begin{equation*} \Pf_{z=0} \Tr([P,Q] \Delta^{-z/r}) = \Res_{z=0} \Tr \left(P \cdot \frac{Q - \Delta^{-z/r} Q \Delta^{z/r}}{z}\Delta^{-z/r} \right) \end{equation*}
by reasoning first for $z \in \C$ of sufficiently large real part to use the trace property, and then applying the analytic continuation property.
\setlength{\parindent}{0mm}
Then, informally we can think of the complex powers of $\Delta$ as 
\begin{equation*} \Delta^{z/r} = e^{\log \Delta \cdot z/r} = 1 + \dfrac{z}{r} \log \Delta + \ldots + \dfrac{1}{p!}\left(\dfrac{z}{r}\right)^p (\log \Delta)^{p} + O(z^{p+1}) \end{equation*}
which after some calculations, gives the expansion 
\begin{equation*}(Q - \Delta^{-z/r} Q \Delta^{z/r})\Delta^{-z/r} = z \delta(Q) - \dfrac{z^2}{2} \delta^2(Q) + \ldots + (-1)^{p-1}\dfrac{z^{p}}{p!}\delta^{p}(Q) + O(z^{p+1}) \end{equation*} 
For a complete proof, see \cite{Rod2013}. But here is a standard example.

\begin{ex} As a more concrete example, let us see how to recover the Noether index theorem from a low dimensional case. Let $M = S^1$ be the unit circle. Consider the operators $D = \frac{1}{\i}\frac{d}{dt}$, $F = D \vert D \vert^{-1}$ and $P = \frac{1+F}{2}$ acting on the Hardy space $H^2(S^1)$. The cosphere bundle of $S^1$ is $S^*S^1 = S^1 \times \{1\} \cup S^1 \times \{-1\}$. Then, remark that $P$ is a pseudodifferential operator of order $0$, its symbol defined on $T^*S^1$ is $\sigma_F(t,\xi) = \frac{1+\xi\vert \xi \vert^{-1}}{2}$, where $\vert . \vert$ denotes the euclidian norm. \\ 
 
Then, let $u \in \cinf(S^1)$ be a nowhere vanishing smooth function. We extend the associated Toeplitz operator $PuP$ to $L^2(S^1)$ by considering the operator $T_u = PuP - (1-P)$, which is an elliptic pseudodifferential of order $0$ of symbol
given by
\[ \left\lbrace \begin{array}{cl} u(t) & \text{ on }S^1 \times (0,\infty) \\ 1 & \text{ on }S^1 \times (-\infty,0) \end{array} \right. \]
Then, using the star-product formula (\ref{star product}), one sees that the part of order $-1$ in the symbol of $T_{u^{-1}}[\log D, T_u]$ is  
\[ \left\lbrace \begin{array}{cl} \dfrac{1}{\i \xi} \dfrac{u'(t)}{u(t)} & \text{ on }S^1 \times (0,\infty) \\ 0 & \text{ on }S^1 \times (-\infty,0) \end{array} \right. \]
Hence:
\[ \Ind(T_u) = - \dfrac{1}{2 \pi \i}\int_{S^1} u^{-1} du  \]
\end{ex}

\section{Relation to the Chern--Connes character}

In this section, we establish the relationship between the Radul cocycle and the Chern--Connes character of a spectral triple.

Let $(\A,H,F)$ be a (trivially graded) $p$-summable Fredholm module. In addition, let $\Psi = \Psi(\Delta)$ be an abstract algebra of pseudodifferential operators, such that 
\begin{enumerate}[label={(\arabic*)}]
\item $\Psi^0$ is an algebra of bounded operators on $H$ containing the representation of $\A$,
\item $\Psi^{-1}$ is a two-sided ideal consisting of $p$-summable operators on $H$, 
\item $F$ is a multiplier of $\Psi^0$ and $[F,\Psi^0]\subset \Psi^{-1}$.
\end{enumerate}

We have an abstract principal symbol exact sequence,  
\begin{equation} \label{ext}
0 \longrightarrow \Psi^{-1} \longrightarrow \Psi^0 \longrightarrow \Psi^0/\Psi^{-1} \longrightarrow 0 
\end{equation}
$\Psi^0/\Psi^{-1}$ should be viewed as an "abstract cosphere bundle". This extension is related to the one involving regularizing operators, as the inclusion of ideals $\psi^{-\infty} \subset \psi^{-1}$ yields the following morphism of extensions,
\begin{equation*}
\begin{tikzpicture}[descr/.style={fill=white}]
\matrix(m)[matrix of math nodes, row sep=3em, column sep=2.5em, text height=2ex, text depth=0.25ex]
    { 0 & \Psi^{-\infty} & \Psi^0 & \S^0=\Psi^0/\Psi^{-\infty} & 0 \\
0 & \Psi^{-1} & \Psi^0 & \Psi^0/\Psi^{-1} & 0 \\};
\draw[double distance = 2pt, line width=0.75pt] (m-1-3) -- (m-2-3) ;
\path[->, font=\scriptsize, >=angle 90] 
(m-1-2) edge (m-2-2)
(m-1-1) edge (m-1-2) 
(m-1-2) edge (m-1-3)
(m-1-3) edge (m-1-4) 
(m-1-4) edge (m-1-5) 
(m-1-4) edge (m-2-4)
(m-1-4) edge (m-1-5)        
(m-2-1) edge (m-2-2)
(m-2-2) edge (m-2-3) 
(m-2-3) edge (m-2-4) 
(m-2-4) edge (m-2-5);        
\end{tikzpicture} 
\end{equation*}

Then, the cyclic cohomology class of the operator trace $[\Tr] \in \HP^0(\Psi^{-\infty})$ extends to a cyclic cohomology class $[\tau] \in \HP^0(\Psi^{-1})$, represented for any choice of integer $k > p$ by the following cyclic $k$-cocycle on $\psi^{-1}$ :
\[ \tau_k(x_0, \ldots, x_k) = \Tr(x_0 \ldots x_k) \]
By naturality of excision, the image of the trace $\partial [Tr] \in \HP^{1}(\S^0)$ by excision is the pull-back of the class $\partial [\tau] \in \HP^{1}(\Psi^0/\Psi^{-1})$. We shall then make a slight abuse of notation by identifying both. \\

Let $P=\frac{1}{2}(1+F)$. Then $[P,a]\in \Psi^{-1}$ for every $a\in A$. The linear map
\begin{equation*}
\rho_F\ :\ \A \longrightarrow \Psi^0/\Psi^{-1}\ ,\qquad \rho_F(a) = PaP \mod \Psi^{-1}\ ,
\end{equation*}
is an algebra homomorphism since $Pa_1Pa_2P = Pa_1a_2 P \mod \Psi^{-1}$ for all $a_1,a_2\in \A$.
\begin{thm} \label{prop Radul Chern}
The Chern-Connes character of the Fredholm module $(H,F)$ is given by the odd cyclic cohomology class over $\A$
\begin{equation*}
\ch(H,F) =\rho^*_F \circ \partial([\Tr]) 
\end{equation*}
where $[\Tr]\in HP^0(\Psi^{-1})$ is the class of the operator trace, $\partial:HP^0(\Psi^{-1}) \to HP^1(\Psi^0/\Psi^{-1})$ is the excision map associated to extension (\ref{ext}), and $\rho_F^*: HP^1(\Psi^0/\Psi^{-1})\to HP^1(\A)$ is induced by the homomorphism $\rho_F$. 
\end{thm}
\begin{pr}
Consider the algebra $\E = \{ (Q,a) \in \Psi^0\oplus \A \, ; \, Q= PaP \mod \Psi^{-1}\}$. The homomorphism $\E \to\A$, $(Q,a)\mapsto a$ yields an extension
\[
0\longrightarrow \Psi^{-1} \longrightarrow \E \longrightarrow \A \longrightarrow 0\ .
\]
The Chern-Connes character $\ch(H,F)\in HP^1(\A)$ is the image of the operator trace by the boundary map of the top this extension. On the other hand, the homomorphism $\E\to\Psi^0$, $(Q,a)\mapsto Q$ yields a commutative diagram of extensions
\[
\begin{tikzpicture}[descr/.style={fill=white}]
\matrix(m)[matrix of math nodes, row sep=3em, column sep=2.5em, text height=2ex, text depth=0.25ex]
    { 0 & \Psi^{-1} & \E & \A & 0 \\
0 & \Psi^{-1} & \Psi^0 & \Psi^0/\Psi^{-1} & 0 \\};
\draw[double distance = 2pt, line width=0.75pt] (m-1-2) -- (m-2-2) ;
\path[->, font=\scriptsize, >=angle 90] 
(m-1-3) edge (m-2-3)
(m-1-1) edge (m-1-2) 
(m-1-2) edge (m-1-3)
(m-1-3) edge (m-1-4) 
(m-1-4) edge (m-1-5) 
(m-1-4) edge node[right]{$\rho_F$} (m-2-4)
(m-1-4) edge (m-1-5)        
(m-2-1) edge (m-2-2)
(m-2-2) edge (m-2-3) 
(m-2-3) edge (m-2-4) 
(m-2-4) edge (m-2-5);        
\end{tikzpicture} 
\]
The conclusion then follows from the naturality of excision. $\hfill{\square}$
\end{pr}

\section{Discussion on manifolds with conical singularities}

Studying index theory on manifolds with singularities is actually one of the motivations for studying a residue index formula adapted to cases where the zeta function exhibits multiple poles. It has indeed been known for many years that zeta functions may exhibit double poles in the context of conical manifolds, see for example the paper of Lescure \cite{Lescure}. In fact, triple poles may also occur, as we shall see. \\

We shall first recall briefly what we need from the theory of conic manifolds, i.e pseudodifferential calculus, residues and results on the associated zeta function. This review part essentially follows the presentation of \cite{GilLoy2002}.

\subsection{Generalities on b-calculus and cone pseudodifferential operators} In our context, manifolds with conical singularities are just manifolds with boundary with an additional structure given by a suitable algebra of differential operators. \\

More precisely, let $M$ be a compact manifold with (connected) boundary, and $r : M \to \R_+$ be a boundary defining function, i.e a smooth function vanishing on $\partial M$ and such that its differential is non-zero on every point of $\partial M$. We work in a collar neighbourhood $[0,1)_r \times \partial M_x$ of the boundary, the subscripts are the notations for local coordinates.  

\begin{df} A \emph{Fuchs type differential operator} $P$ of order m is a differential operator on $M$ which can be written in the form
\begin{equation*} P(r,x) = r^{-m} \sum_{j+\vert \alpha \vert \leq m} a_{j,\alpha}(r,x) (r \partial_r)^j \partial_x^\alpha \end{equation*} 
in the collar $[0,1)_r \times \partial M_x$. The space of such operators will be denoted $r^{-m} \Diff_b^m(M)$. 
\end{df}

$\Diff_b^m(M)$ denotes the algebra of $b$-differential operators in Melrose's calculus for manifolds with boundary. We now recall the associated small $b$-pseudodifferential calculus $\Psi_b(M)$.\\

Let $M_b^2$ be the $b$-stretched product of $M$, i.e the manifold with corners whose local charts are given by the usual charts on $M^2 \smallsetminus \partial M^2$, and parametrized by polar coordinates over $\partial M$ in $M^2$. More precisely, writing $M \times M$ near $r=r'=0$ as
\[ M^2 \simeq [0,1]_r \times [0,1]_{r'} \times \partial M^2 \]
this means that we parametrize the part $[0,1]_r \times [0,1]_{r'}$ in polar coordinates
\[ r = \rho \cos \theta , \quad r' = \rho \sin \theta \] 
for $\rho \in \R_+$, $\theta \in [0,\pi/2]$. The right and left boundary faces are respectively the points where $\theta = 0$ and $\theta = \pi/2$. \\

Let $\Delta_b$ the $b$-diagonal of $M^2_b$, that is, the lift of the diagonal in $M^2$. Note that $\Delta_b$ is in fact diffeomorphic to $M$, so that any local chart on $\Delta_b$ can be considered as a local chart on $M$.  \\ 
%[make a picture ...] \\
\begin{df} The \emph{algebra of $b$-pseudodifferential operators of order $ m $}, denoted $\Psi_b^m(M)$, consists of operators $P : C^\infty(M) \to C^\infty(M)$ having a Schwartz kernel $K_P$ such that 
\begin{enumerate}[label={(\roman*)}]
\item Away from $\Delta_b$, $K_P$ is a smooth kernel, vanishing to infinite order on the right and left boundary faces. \\

\item On any local chart of $M^2_b$ intersecting $\Delta_b$ of the form $U_{r,x} \times \R^n$ such that $\Delta_b \simeq U \times \{0\}$, and where $U$ is a local chart in the collar neighbourhood $[0,1)_r \times \partial M_x$ of $\partial M$, we have
\begin{equation*}
K_P(r,x,r',x') = \dfrac{1}{(2\pi)^n} \int e^{\i (\log(r/r') \cdot \tau + x \cdot \xi)} a(r,x,\tau,\xi) \, d\tau \, d\xi
\end{equation*} 
where $a(y,\nu)$, with $y = (r,x)$ and $\nu = (\tau, \xi)$, is a classical pseudodifferential symbol of order $m$, plus the condition that $a$ is smooth in the neighbourhood of $r=0$. 
\end{enumerate}
\end{df}

Remark that $\log(r/r')$ should be singular at $r = r' = 0$ if we would have considered kernels defined on $M^2$. Introducing the $b$-stretch product $M^2_b$ has the effect of blowing-up this singularity. \\

The \emph{algebra of conic pseudodifferential operators} is then the algebra $r^{-\Z} \Psi_b^{\Z}(M)$. The opposed signs in the filtrations are only to emphasize that $r^{\infty}\Psi_b^{-\infty}(M)$ is the associated ideal of regularizing operators. \\

To such an operator $A = r^{-p} P \in r^{-p} \Psi_b^{m}$ , we define on the chart $U$ the local density 
\begin{equation*}
\omega(P)(r,x) = \left(\int_{\vert \nu \vert = 1} p_{-n}(r,x,\tau,\xi) \iota_L d\tau d\xi \right) \cdot \, \dfrac{dr}{r} \, dx
\end{equation*}  
where $\nu = (\tau, \xi)$ and $L$ is the generator of the dilations. \\

It turns out (but this is not obvious) that this a priori local quantity does not depend on the choice of coordinates on $M$, and hence, define a globally defined density $\omega(P)$, smooth on $M$, that we call the \emph{Wodzicki residue density}. Unfortunately, the integral on $M$ of this density does not converge in general, as the boundary introduces a term in $1/r$ in the density. However, we can regularize this integral, thanks to the following lemma. Here, $\Omega_b$ denote the bundle of $b$-densities on $M$, that is, the trivial line bundle with local basis on the form $(dr/r)dx$. The following lemma from Gil and Loya is proved in \cite{GilLoy2002}.   
 
\begin{lem} \label{b-regularization} Let $r^{-p} u \in C^\infty(M,\Omega_b)$, and $p \in \R$. Then, the function
\begin{equation*}
z \in \C \longmapsto \int_M r^z u
\end{equation*}
is holomorphic on the half plane $\Re z > p$, and extends to a meromorphic function with only simple poles at $z=p, p-1, \ldots $. If $p \in \N$, Its residue at $z=0$ is given by
\begin{equation}
\Res_{z=0} \int_M r^z u(r,x) \, \frac{dr}{r} dx = \dfrac{1}{p!} \int_{\partial M} \partial_r^p (r^p u(r,x))_{r=0} \, dx 
\end{equation}
\end{lem}

Applying this regularization to the Wodzicki residue density is useful to many "residues traces" that we immediately study. 

\subsection*{Traces on conic pseudodifferential operators}

We first begin by defining different algebras of pseudodifferential operators, introduced by Melrose and Nistor in \cite{MelNis1996}. The main algebra that we shall consider is 
\begin{equation*}
A = r^{-\Z} \Psi^\Z(M) = \bigcup_{p \in \Z} \bigcup_{m \in \Z} r^{- p} \Psi^m(M)
\end{equation*}
which clearly contains the algebra of Fuchs type operators. The ideal of \emph{regularizing operators} is
\begin{equation*}
I = r^{\infty} \Psi^{-\infty}(M) = \bigcup_{p \in \Z} \bigcup_{m \in \Z} r^{- p} \Psi^m(M)
\end{equation*} 
and this explains why we note the two filtrations by opposite signs in $A$. Consider the following quotients
\begin{equation*}
I_{\sigma} = r^{\infty} \Psi^\Z(M)/I, \quad I_{\partial} = r^\Z \Psi^{-\infty}(M)/I
\end{equation*} 
Here, $I_{\sigma}$ should be thought as an extension of the algebra of pseudodifferential operators in the interior of $M$, whereas $I_{\partial}$ are     
regularizing operators up to the boundary. We finally define 
\begin{equation*}
A_\partial = \A / I_{\sigma}, \quad A_\sigma = A / I_{\partial}, \quad A_{\partial,\sigma} = A / (I_{\partial} + I_{\sigma})
\end{equation*}

\begin{df} Let $P \in r^{-p} \Psi^m(M)$ be a conic pseudodifferential operator, with $p,m \in \Z$. According to Lemma \ref{b-regularization}, define the functionals $\Tr_{\partial, \sigma}$, $\Tr_{\sigma}$ to be
\begin{gather}
\Tr_{\partial, \sigma}(P) = \Res_{z=0} \int_M r^z \omega(P)(r,x) \, \dfrac{dr}{r} dx = \dfrac{1}{p!} \int_{\partial M} \partial_r^p (r^p \omega(P)(r,x))_{r=0} \, dx \\
\Tr_{\sigma}(P) = \Pf_{z=0} \int_M r^z \omega(P) \, \dfrac{dr}{r} dx
\end{gather}
where $\Pf$ denotes the constant term in the Laurent expansion of a meromorphic function.  
\end{df}

\begin{rk} Using Lemma \ref{b-regularization}, one can show that $\Tr_{\partial, \sigma}(P)$ does not depend on the choice of the boundary defining function $r$. This is not the case for $\Tr_{\sigma}(P)$, but its dependence on $r$ can be explicitly determined, cf. \cite{GilLoy2002}. 
\end{rk}

The "Partie Finie" regularization of a trace does not give in general a trace, and this is indeed the same for the functional $\Tr_{\sigma}(P)$ acting on these algebras, the obstruction to that is precisely the presence of the boundary. However, by definition, $\Tr_{\sigma}(P)$ clearly defines an extension of the Wodzicki residue for pseudodifferential operators, one can expect that it is a trace on $I_{\sigma} = r^{\infty} \Psi^\Z(M)/I$. 

\begin{thm} \emph{(Melrose-Nistor, \cite{GilLoy2002, MelNis1996})} $\Tr_{\sigma}$ is, up to a multiplicative constant, the unique trace on the algebra $I_\sigma$
\end{thm}  

By Lemma \ref{b-regularization} and the definition above, the defect of $\Tr_{\sigma}$ to be a trace is precisely measured by $\Tr_{\partial, \sigma}(P)$, which can therefore be viewed as a restriction of the Wodzicki residue to the boundary $\partial M$. Then, the following proposition seems natural.  

\begin{thm} \emph{(Melrose-Nistor, \cite{GilLoy2002, MelNis1996})} $\Tr_{\partial, \sigma}$ is, up to a multiplicative constant, the unique trace on the algebras $A_\partial$, $A_\sigma$ and $A_{\partial,\sigma}$
\end{thm}

These two traces may be seen as "local" terms, since they only depend on the symbol of the pseudodifferential operator considered. The first can be seen as a trace on interior of $M$, the second is related to the boundary $\partial M$. There is one last trace to introduce, less easy to deal with because this one is not local. \\

Fix a holomorphic family $Q(z) \in r^{\alpha z} \Psi_b^{\beta z}(M)$, with $\alpha, \beta \in \R$, such that $Q$ is the identity at $z=0$. Take $P \in r^{-p} \Psi_b^m$, with $p,m \in \Z$ and let $(PQ(z))_\Delta$ be the restriction to the diagonal $\Delta$ of $M^2$ of the Schwartz kernel of $PQ(z)$. Melrose and Nistor noticed in \cite{MelNis1996} that $(PQ(z))_\Delta$ is meromorphic in $\C$, with values in $r^{\alpha z - p} C^\infty(M)$ with possible simple poles in the set 
\[ \left\{\dfrac{-n-m}{\beta}, \dfrac{-n-m+1}{\beta}, \ldots \right\}  \]

\begin{df} Let $P \in r^{-p} \Psi_b^m$ be a conic pseudodifferential operator. Then, we define 
\begin{equation*}
\Tr_\partial(P) = \dfrac{1}{p!} \int_{\partial M} \partial_r^p (r^p \Pf_{z=0} (PQ(z))_\Delta)_{r=0} \, dx 
\end{equation*}
If $p$ is not an integer, then, $\Tr_\partial(P)$ is defined to be $0$. 
\end{df}

\begin{rk} $\Tr_\partial(P)$ depend on the choice of the operator $Q$, but the dependence can be explicitly determined, see \cite{MelNis1996}.
\end{rk}

There is an interpretation of  $\Tr_\partial$ analogous to those of $\Tr_{\partial, \sigma}$ : If the order of $P$ is less than the dimension of $M$, then $\Tr_\partial(P)$ is a kind of $L^2$ of $P$ restricted to the boundary. This is precisely the content of the following result. 

\begin{thm} \emph{(Melrose-Nistor, \cite{GilLoy2002, MelNis1996})} $\Tr_\partial(P)$ is, up to a multiplicative constant, the unique trace on the algebra \[ I_{\partial} = r^\Z \Psi^{-\infty}(M)/I \]
\end{thm}

\subsection*{Heat kernel expansion and zeta function}

Now, let  $\Delta \in r^{-2} \Diff_b^{2}(M)$ be \emph{fully elliptic}, or \emph{parameter elliptic} with respect to a parameter $\alpha$. We refer to \cite{GilLoy2002} for the definition, what we need to know is just that this condition ensures the existence of the heat kernel $e^{-t\Delta}$ of $A$, and that operators of the type $P \Delta^{-z}$, with $P \in r^{-p} \Psi_b^m$, are of trace-class on $r^{\alpha - m}L^2_b(M)$ for $z$ in the half-plane $\Re z > max\{\frac{m + n}{2} , \frac{p}{2} \}$, $n = \dim \, M$. 

\begin{ex} \label{conic laplacian} As usual, we work in a collar neighbourhood of $M$. Then, the operator 
\begin{equation} \label{ex fully elliptic}  
\Delta = \frac{1}{r^2} \left( (r \partial_r)^2 - \Delta_{\partial M} + \frac{(n-2)^2}{4} + a^2 \right) 
\end{equation} 
where $a > 1$, is and $\alpha = 1$, is an example of such an operator. See \cite{GilLoy2002} for more details.
\end{ex}

Then, the traces introduced in the previous paragraph gives the coefficients of the expansion of $\Tr(P e^{-t\Delta})$. 

\begin{thm} \emph{(Gil-Loya, \cite{GilLoy2002})}  Under the conditions above, we have 
\begin{equation*} 
\Tr(P e^{-t\Delta}) \sim_{t \to 0} \sum_{k \geq 0} a_k t^{(k-p)/2} + (b_k + \beta_k \log t)t^{k} + (c_k + \gamma_k \log t + \delta_k (\log t)^2)t^{(k-m-n)/2}
\end{equation*}   
where 
\begin{align*}
&\beta_k = C_k (\Tr_\sigma + \Tr_\partial)(P \Delta^k) \\
&\gamma_k = C'_K \Tr_{\partial, \sigma} (P \Delta^{k-m-n}) \\
&\delta_k = C_k'' \Tr_{\partial, \sigma} (P \Delta^{k-m-n})
\end{align*}
$C_k$, $C'_K$, $C_k''$ are explicit (but not of interest for us). \\

In particular, the coefficient of $\log t$ is 
\[ -\dfrac{1}{2} \Tr_\sigma(P) - \dfrac{1}{2} \Tr_\partial(P) - \dfrac{1}{4} \Tr_{\partial, \sigma}(P) \]
and the coefficient of $(\log t)^2$ is
\[ - \dfrac{1}{4} \Tr_{\partial, \sigma}(P) \]
\end{thm}

Using a Mellin transform, we can write
\[ \Tr(P \Delta^{-z/2} = \dfrac{1}{\Gamma(z/2)} \int_0^{\infty} t^{z-1} \Tr(P e^{-t\Delta}) \, dt \]
and knowing, that $z \mapsto \int_1^{\infty} t^{z-1} \Tr(P e^{-t\Delta}) \, dt$ is entire, the asymptotic expansion of the previous proposition gives the following corollary on the zeta function. 

\begin{cor} \label{conic zeta} The zeta function $z \mapsto \Tr(P \Delta^{-z/2})$ is holomorphic in the half-plane $\Re z > max\{m + n , p \}$, and extends to a meromorphic function with at most triple poles, whose set is discrete. At $z=0$, there are simple and double poles only, which are respectively given by the terms of $\log t$ and $(\log t)^2$ in the heat kernel expansion of $\Tr(P e^{-t\Delta})$. 
\end{cor}
 
\subsection{Spectral triple and regularity} In this paragraph, we want to investigate if Fuchs type operators on conic manifolds can define an abstract algebra of differential operators, so that the local index formula we gave in the first section applies. \\

We start with a conic manifold. Let $M$ be a manifold with connected boundary, with boundary defining function $r$, endowed with the algebra of Fuchs type differential operators. The points (i), (ii), (iii) of Definition \ref{df ADO} are verified, if for example we take for $\Delta$ the fully-elliptic operator of order 2 given in Example \ref{ex fully elliptic}, and require that the order is given by the differential order. More generally, working locally in a collar neighbourhood $[0,1)_r \times \partial M_x$ of the boundary $\partial M$, elementary calculations shows that
\begin{equation} \label{commutator Fuchs order}
[r^p \Diff^m_b(M), r^{p'} \Diff^{m'}_b(M)] \subset r^{p+p'} \Diff^{m + m' -1}_b(M) 
\end{equation}
and as we shall see, the fact that the order in $r$ does not decrease is the problem. \\

Let us denote by $r^p C^\infty(\partial M)$ the subalgebra of $C^\infty(M)$ of functions $f$ which have an asymptotic expansion 
\[ f(r,x) \sim r^p f_p(x) + r^{p+1} f_{p+1}(x) + \ldots \] 
in a neighbourhood of $r=0$. Here, the $\sim$ means that the rest of such an expansion is of the form $r^N f_N(r,x)$, with $f_N$ bounded in the collar $[0,1) \times \partial M$. The case $p=0$ actually corresponds to the smooth functions on the collar. \\ 

For the algebra of the spectral triple, it seems a good choice to look for a candidate among these classes of functions. But doing so, the formula of Lemma \ref{CM trick} is no more asymptotic in the sense of Definition \ref{asymptotic expansion}. Indeed, if $b(r,x) = r^p$ for $p \in \N$, the observation (\ref{commutator Fuchs order}) shows that the terms $b^{(k)}$ are in $r^{p-2k} \Diff_b^{k}(M)$, but by the properties of the zeta function given in the Corollary \ref{conic zeta}, the function
\[ z \longmapsto \Tr(b^{(k)}\Delta^{-k-z}) \]
is holomorphic for $\Re(z) + k > \max\left\{\frac{n+k}{2}, \frac{2k - p}{2} \right\}$, which is equivalent to $\Re(z) > \max\left\{\frac{n-k}{2}, -\frac{p}{2} \right\}$. Hence, if $p \geq 0$, the function above is in general not holomorphic at $0$ when $N$ goes to infinity. In other terms, the spectral triple we may construct will be not regular, and local index formulas of Connes-Moscovici, or those given at the beginning cannot be applied directly. As we have seen, the main problem is due to the fact that there are two notions of order : The differential order, which is local, and "the order in $r$", which is not, and comes form the presence of the boundary $\partial M$. \\

However, we may recover some interesting informations on $M$ from the zeta function. Note for instance that the higher residue $\littlebarint^2$ defined in Proposition \ref{higher WG trace} gives the trace $\Tr_{\partial, \sigma}$. $\littlebarint^1$ is, modulo some constant terms, the sum of the three functionals $\Tr_{\partial, \sigma}$, $\Tr_{\sigma}$, $\Tr_{\partial}$, which illustrates that it is no more a trace on the algebra of conic pseudodifferential operators. The next paragraph is a discussion on index theory. 

\subsection{A non-local index formula}

The formula of Theorem \ref{local index formula} cannot be applied directly since we are not in the context of regular spectral triples. However, there are always some relevant informations to get on index theory. \\

Let $M$ be a manifold with boundary, seen as a conic manifold, and consider the extension 
\[ 0 \to r^{\infty} \Psi_b^{-\infty}(M) \to r^{-\Z} \Psi_b^{\Z}(M) \to r^{-\Z} \Psi_b^{\Z}(M)/r^{\infty} \Psi_b^{-\infty}(M) \to 0 \] 
Here, by an \emph{elliptic pseudodifferential operator} $P \in r^{-\Z} \Psi_b^{\Z}(M)$, we shall mean that $P$ is invertible in the quotient $A = r^{-\Z} \Psi_b^{\Z}(M)/r^{\infty} \Psi_b^{-\infty}(M)$. Being \emph{fully elliptic} is an extra condition on the {indicial or normal operator}, which guarantees that $P$ is Fredholm between suitable spaces. We shall not enter into these details : What we want to investigate is just the pairing given in the paragraph (\ref{index pairing}). In particular, if $P$ is fully elliptic, then the pairing really calculates a Fredholm index. \\
 
Now, let $P,Q \in r^{-\Z} \Psi_b^{\Z}(M)$.  We can still follow the "Partie Finie" argument given in the proof of Theorem \ref{local index formula}, so that we still have the Radul cocycle
\begin{align*} 
c(P,Q) & = \Pf_{z=0} \Tr([P,Q] \Delta^{-z})  \\
& \Res_{z=0} \Tr \left(P \cdot \left( \dfrac{Q - \Delta^{-z} Q \Delta^{-z}}{z}\right)\Delta^{-z}\right) 
\end{align*} 
As we already said, the Connes-Moscovici's formula in Lemma \ref{CM trick} is no more asymptotic, but from an algebraic viewpoint, the (\ref{CM trick bis}) still holds. So, for any integer $N$, which will be thought large enough, we have  
\begin{equation*}
Q - \Delta^{-z} Q \Delta^{-z} = \sum_{k=1}^N Q^{(k)} \Delta^{-k} + \dfrac{1}{2\pi\i} \int \lambda^{-z} (\lambda - \Delta)^{-1} Q^{(N+1)} (\lambda - \Delta)^{-N-1} \, d\lambda 
\end{equation*} 
We now take advantage of the fact that the traces $\Tr_{\sigma}$ and $\Tr_{\partial, \sigma}$ vanishes when the differential order of the operators is less that the dimension of $M$. We then have the following result.

\begin{thm} Let $M$ be a conic manifold, i.e a manifold with boundary endowed with a conic metric, and let $r$ be a boundary defining function. Let $\Delta$ be the "conic laplacian" of Example \ref{conic laplacian}. Then, the Radul cocycle associated to the pseudodifferential extension  
\[ 0 \to r^{\infty} \Psi_b^{-\infty}(M) \to r^{-\Z} \Psi_b^{\Z}(M) \to r^{-\Z} \Psi_b^{\Z}(M)/r^{\infty} \Psi_b^{-\infty}(M) \to 0 \]
is given by the following \emph{non local} formula :
\begin{multline*}
c(a_0,a_1) = (\Tr_{\partial, \sigma} + \Tr_{\sigma})(a_0[\log \Delta, a_1]) - \frac{1}{2}\Tr_{\partial, \sigma}(a_0[\log \Delta,[\log \Delta, a_1]]) +  \\
+ \Tr_\partial \left( a_0 \sum_{k=1}^N a_1^{(k)} \Delta^{-k} \right) + \dfrac{1}{2\pi\i} \Tr\left(\int \lambda^{-z} a_0 (\lambda - \Delta)^{-1} a_1^{(N+1)} (\lambda - \Delta)^{-N-1}\right) \, d\lambda
\end{multline*} 
for $a_0, a_1 \in \Psi_b^{\Z}(M)/r^{\infty} \Psi_b^{-\infty}(M)$
\end{thm}
In the right hand-side, the first line consists in local terms only depending on the symbol of $P$, the second line gives the non local contributions. \\

If $P \in r^{-\Z} \Psi_b^{\Z}(M)$ is an elliptic operator, so that $P$ defines an element in the odd K-theory group $K_1^{\alg}(A)$, and $Q$ an inverse of $P$ modulo $A$, we then obtain a formula for the index of $P$. The second line of the formula above should be a part of the eta invariant (when it is defined). A perspective may be to investigate how to compare these different elements in order to get another definition of the eta invariant, suitable not only for Dirac operators but also for general pseudodifferential operators.   

\bibliographystyle{plain}
\bibliography{bibliographie}

\end{document}